# Summation functions with nonlinear asymptotic behavior

Victor Volfson

ABSTRACT The paper considers a universal approach that allows one to quite simply obtain nonlinear asymptotic estimates of various summation functions. It is shown the application of this approach to the asymptotic estimation of divergent Dirichlet series. Several assertions have been proven and numerous examples have been considered.

Key words: arithmetic function, summation function, linear asymptotic of the summation function, nonlinear asymptotics of the summation functions, Abel summation formula, finite limit of the average value of an arithmetic function, asymptotics of divergent Dirichlet series.

# 1. INTRODUCTION

In general, an arithmetic function $f(m), m = 1,...,n$ is a function defined on the set of natural numbers and taking values on the set of complex numbers. The name arithmetic function is due to the fact that this function expresses some arithmetic property of the natural series.

A function of the form:

$$S(n) = \sum_{m \leq n} f(m). \tag{1.1}$$

is called summation.

The Mertens function is a summation function - $M(n) = \sum_{m \leq n} \mu(m)$, where $\mu(m)$ is the arithmetic Möbius function.

The assertion was proven in [1]. Let there is a finite limit

$$d^*(S) = \lim_{n \to \infty} \frac{S(n)}{n}, \tag{1.2}$$

where $d^*(S)$ is the asymptotic density jf a summation function $S(n)$. Then the asymptotic behavior of the summation function $S(n)$ at $n \to \infty$ has the form:

$$S(n) = d^*(S)n + o(n). \tag{1.3}$$

Let us call asymptotic (1.3) linear asymptotic.

Consequently, the existence of a finite limit (1.2) is necessary (for the summation function) to have linear asymptotic behavior (1.3).

Please note that based on (1.1), limit (1.2) can be written as:

$$d^*(S) = \lim_{n \to \infty} \frac{\sum_{m \leq n} f(m)}{n}, \tag{1.4}$$

where $\frac{1}{n} \sum_{m \leq n} f(m)$ is the average value of the arithmetic function $f(m), m = 1,...,n$.

Based on (1.4), we can formulate that the existence of a finite limit of the average value of the arithmetic function $f(m), m = 1,...,n$ at $n \to \infty$ is necessary in order for the summation function to have linear asymptotic behavior (1.3).

Summation functions with a finite mean limit (respectively, linear asymptotic behavior (1.3)) were studied in [2].

We will consider summation functions whose asymptotic estimates are nonlinear in this work. For example, the following asymptotic equality is true for the summation functions of the Euler arithmetic function $\sum_{n \leq x} \varphi(n) \sim \frac{3}{\pi^2} x^2, x \to \infty$.

Typically, methods for obtaining such estimates are quite difficult [3,4]. Therefore, our goal is to obtain an approach (method) that would make it possible to relatively simply obtain such estimates for a fairly large class of arithmetic functions.

Currently, unfortunately, there is no such approach. This explains the various and difficult ways of obtaining these estimates in works on this topic. Therefore, this approach is original, has scientific novelty, and its application will allow us to approach the solution of this problem from a unified scientific position.

A detailed and meaningful consideration of various cases in the work will make it possible to obtain the indicated asymptotic estimates for various summation functions using fairly simple methods.

A large number of examples show the practical application of this approach to obtain asymptotics of various summation functions.

2. FINDING NONLINEAR ASYMPTOTIC BEHAVIOR OF SUMMATION FUNCTIONS

The approach implemented in this work is based on average estimates.

Suppose that some arithmetic function $f(n)$ is asymptotically equivalent to a monotonic arithmetic function, i.e. $f(n) \sim g(n)$. Then the following asymptotic estimate on average for the summing function is true:

$$\sum_{n \leq x} f(n) \sim \sum_{n \leq x} g(n), x \to \infty. \tag{2.1}$$

Another tool that will be used in this work is Abel's summation formula:

$$\sum_{n \leq x} a_n b(n) = A(x)b(x) - \int_1^x A(t)b'(t)dt, \tag{2.2}$$

where $a_n$ is some sequence of real or complex numbers, $A(x) = \sum_{n \leq x} a_n$, $b(t)$ is some continuously differentiable function on [1, x].

This formula is a discrete analogue of the method of integration over parts.

The next tool that will be used in the work is the formula for the linear asymptotic of the summation function (1.3).

Assertion1

Let there is an arithmetic function $f(n), n = 1,..., x$ for which there is a finite limit of the average value:

$$\lim_{n \to \infty} \frac{1}{n} \sum_{n \leq x} f(n) = c, \tag{2.3}$$

where $c \neq 0$, then the following asymptotic equality holds:

$$\sum_{n \leq x} nf(n) \sim \frac{c}{2} x^2, x \to \infty. \tag{2.4}$$

Proof

Let us take the arithmetic function $f(n), n = 1,..., x$ as a sequence $a_n$, then in accordance with (1.3) $A(x) = cx + o(x)$. Let us put $b(n) = n$. We substitute everything into (2.2) and get:

$$\sum_{n \leq x} nf(n) = (cx + o(x))x - \int_1^x (ct + o(t))dt = cx^2 + o(x^2) - \frac{cx^2}{2} + O(1) + o(x^2) = \frac{c}{2}x^2 + o(x^2)$$

which corresponds to (2.4).

As an example of Assertion 1, we consider finding the asymptotic of the summation function $\sum_{n \leq x} \varphi(n)$.

Let's use the fact that $\varphi(n) = n\prod_{p|n}(1-1/p)$. In this case, $f(n) = \prod_{p|n}(1-1/p)$ is a multiplicative function, and $f(n) = \prod_{p|n}(1-1/p) \leq 1$. Taking this into account, based on Wirsing's theorem [5]:

$$c = M(f) = \prod_p (1-\frac{1}{p})(1+\sum_{\alpha \geq 1}\frac{f(p^\alpha)}{p^\alpha}) = \prod_p(1-\frac{1}{p^2}) = \frac{6}{\pi^2} \neq 0.$$

Therefore, based on (2.3) and (2.4):

$$\sum_{n \leq x}\varphi(n) \sim \frac{c}{2}x^2 = \frac{3}{\pi^2}x^2, x \to \infty. \qquad (2.5)$$

Now let's consider another example for assertion 1 - finding the asymptotic of the summation function $\sum_{n \leq x}(\Omega(n) - \omega(n))n$. $f(n) = \Omega(n) - \omega(n)$ is an additive arithmetic function in this case.

Let's define the value:

$$c = \lim_{x \to \infty}\frac{1}{x}\sum_{n \leq x}(\Omega(n) - \omega(n)) = \lim_{x \to \infty}\frac{1}{x}\sum_p\sum_{v \geq 2}\left\lfloor\frac{x}{p^v}\right\rfloor = \sum_p\sum_{v \geq 2}\frac{1}{p^v} = \sum_p\frac{1}{p(p-1)} \neq 0$$

Therefore, based on (2.3) and (2.4):

$$\sum_{n \leq x}(\Omega(n) - \omega(n))n \sim \frac{c}{2}x^2 = \frac{x^2}{2}\sum_p\frac{1}{p(p-1)}, x \to \infty.$$

Assertion 2

Let there is an arithmetic function $f(n), n = 1,...,x$ for which the linear asymptotic behavior with an estimate of the remainder term is known:

$$\sum_{n \leq x}f(n) = cx + O(g(x)), x \to \infty, \qquad (2.6)$$

where $c \neq 0$, then the following asymptotic is true:

$$\sum_{n \leq x}nf(n) = \frac{cx^2}{2} + O(xg(x)) + O(\int_1^x g(t)dt) + O(1), x \to \infty. \qquad (2.7)$$

Proof

Let us put $A(x) = cx + O(g(x))$, $b(n) = n$ in formula (2.2) and get:

$$\sum_{n \leq x} nf(n) = (cx + O(g(x)))x - \int_1^x (ct + O(g(t)))dt = cx^2 + O(xg(x)) - \frac{cx^2}{2} + \frac{c}{2} + O(\int_1^x g(t)dt) = \frac{cx^2}{2} + O(xg(x)) + O(\int_1^x g(t)dt) + O(1)$$

which corresponds to (2.7).

Let us consider examples of obtaining asymptotic estimates based on Assertion 2.

The linear asymptotic behavior of the summation function is known: $\sum_{n \leq x} \frac{\sigma(n)}{n} = \zeta(2)x + O(\ln x), x \to \infty$. This asymptotic satisfies condition (2.6), therefore, to obtain the nonlinear asymptotic of the summation function we use (2.7):

$$\sum_{n \leq x} \sigma(n) = \frac{\zeta(2)x^2}{2} + O(x \ln x) + O(\int_1^x \ln(t)dt) + O(1) = \frac{\pi^2}{12}x^2 + O(x \ln x) + O(x \ln x) + O(1) = \frac{\pi^2}{12}x^2 + O(x \ln x)$$

at $x \to \infty$.

Another example. Let $K(n) = \prod_{p|n} p$ and the linear asymptotic of the summation function is known $\sum_{n \leq x} \frac{K(n)}{n} = Cx + O(\sqrt{x})$, $x \to \infty$, where $C = \prod_p (1 - \frac{1}{p(p+1)}) \neq 0$. This asymptotic behavior satisfies condition (2.6), therefore, to obtain a nonlinear asymptotic of the summation function we use (2.7):

$$\sum_{n \leq x} K(n) = \frac{Cx^2}{2} + O(x\sqrt{x}) + O(\int_1^x \sqrt{t}dt) + O(1) = \frac{C}{2}x^2 + O(x\sqrt{x}) + O(x\sqrt{x}) + O(1) = \frac{C}{2}x^2 + O(x\sqrt{x})$$

at $x \to \infty$.

Formula (2.4) can be generalized.

Assertion 3

Let there is an arithmetic function $f(n), n = 1, ..., x$ for which there is final limit of average value:

$$\lim_{x \to \infty} \frac{1}{x} \sum_{n \leq x} f(n) = c, \qquad (2.8)$$

where $c \neq 0$, then the asymptotic equality is true:

$$\sum_{n \leq x} n^s f(n) \sim \frac{c}{s+1} x^{s+1}, x \to \infty, \qquad (2.9)$$

where $s \geq 1$.

Proof

Based on (2.6), we put $A(x) = cx + o(x)$, $b(n) = n^s, s \geq 1$ in formula (2.2) and get:

$$\sum_{n \leq x} n^s f(n) = (cx + o(x))x^s - s\int_1^x (ct + o(t))t^{s-1} dt = cx^{s+1} + o(x^{s+1}) - \frac{sc}{s+1} x^{s+1} + o(x^{s+1}) \sim \frac{c}{s+1} x^{s+1},$$

for $x \to \infty$, which corresponds to (2.9).

Let's consider an example for assertion 3. Let's determine the asymptotic of the summation function $\sum_{n \leq x} Q(n) n^s = \sum_{n \leq x} \frac{Q(n)}{n} n^{s+1}$.

It is known that $Q(n) = \frac{6}{\pi^2} n + o(n)$, therefore:

$$c = \lim_{x \to \infty} \frac{1}{x} \sum_{n \leq x} \frac{Q(n)}{n} = \lim_{x \to \infty} \frac{1}{x} \sum_{n \leq x} (\frac{6}{\pi^2} + o(1)) = \frac{6}{\pi^2} \neq 0.$$

Since the condition (2.8) is satisfied, then using (2.9), we obtain:

$$\sum_{n \leq x} Q(n) n^s = \sum_{n \leq x} \frac{Q(n)}{n} n^{s+1} \sim \frac{cx^{s+2}}{s+2} = \frac{6x^{s+2}}{\pi^2 (s+2)}$$

at $x \to \infty$.

If the nonlinear asymptotic of a summation function with a remainder term is known, then it can also be used to obtain more detailed asymptotic of another summation function.

Assertion 4

Let there is an arithmetic function $f(n), n = 1, ..., x$ for which nonlinear asymptotic behavior with a remainder term is known:

$$\sum_{n \leq x} f(n) = g_1(x) + O(g_2(x)), \qquad (2.10)$$

for $x \to \infty$, then the following asymptotic behavior of the summation function is true:

$$\sum_{n \leq x} nf(n) = xg_1(x) + O(xg_2(x)) - \int_1^x g_1(t)dt + O(\int_1^x g_2(t)dt) \qquad (2.11)$$

at $x \to \infty$.

Proof

Let us put the value $A(x) = g_1(x) + O(g_2(x))$ and $b(n) = n$ in formula (2.2) in this case and we get:

$$\sum_{n \leq x} nf(n) = (g_1(x) + O(g_2(x)))x - \int_1^x (g_1(t)dt + O(g_2(t)))dt = xg_1(x) - \int_1^x (g_1(t)dt + O(xg_2(x)) + O(\int_1^x g_2(t)dt),$$

for $x \to \infty$, which corresponds to (2.11).

Let's consider an example for assertion 4. It is required to find the asymptotic of the summation function $\sum_{n \leq x} n2^{\omega(n)}$.

The nonlinear asymptotic behavior of the summation function with a remainder term is known $\sum_{n \leq x} 2^{\omega(n)} = \frac{6}{\pi^2} x \ln x + O(x)$ at $x \to \infty$, which corresponds to (2.10). Since based on (2.11) we obtain:

$$\sum_{n \leq x} n2^{\omega(n)} = \frac{6}{\pi^2} x^2 \ln x + O(x^2) - \frac{6}{\pi^2} \int_1^x t \ln t \, dt + O(\int_1^x t \, dt) = \frac{3x^2}{\pi^2} \ln x + O(x^2)$$

at $x \to \infty$.

Let us generalize the assertion 4.

Assertion 5

Let there is an arithmetic function $f(n), n = 1, ..., x$ for which we know nonlinear asymptotic with remainder term:

$$\sum_{n \leq x} f(n) = g_1(x) + O(g_2(x)), \qquad (2.12)$$

at $x \to \infty$, then the following asymptotic of the summation function is true (at $s \geq 1$):

$$\sum_{n \leq x} n^s f(n) = x^s g_1(x) + O(x^s g_2(x)) - s\int_1^x g_1(t)dt + O(\int_1^x g_2(t)t^{s-1}dt). \qquad (2.13)$$

at $x \to \infty$.

Proof

Based on (2.12), we put the value $A(x) = g_1(x) + O(g_2(x))$ and $b(n) = n^s$ in formula (2.2), and we get:

$$\sum_{n \leq x} n^s f(n) = x^s g_1(x) + O(x^s g_2(x)) - s\int_1^x g_1(t)dt + O(\int_1^x g_2(t)t^{s-1}dt) =$$

$$= x^s g_1(x) + O(x^s g_2(x)) - s\int_1^x g_1(t)dt + O(\int_1^x g_2(t)t^{s-1}dt)$$

for $x \to \infty$, which corresponds to (2.13).

Let's consider an example for assertion 5. We need to find the asymptotic behavior of the summation function $\sum_{n \leq x} n^s \tau(n)$.

Having in mind that the nonlinear asymptotic behavior of the summation function $\sum_{n \leq x} \tau(n) = x \ln x + O(x)$ at $x \to \infty$ corresponds to (2.12), then based on (2.13), we obtain:

$$\sum_{n \leq x} n^s \tau(n) = x^{s+1} \ln x + O(x^{s+1}) - s\int_1^x t^s \ln t\, dt + O(\int_1^x t^s dt) = \frac{s}{s+1} x^{s+1} \ln x + O(x^{s+1})$$

at $x \to \infty$.

3. APPLICATION OF THIS APPROACH TO FINDING ASYMPTOTICS OF DIVERGENT DIRICHLET SERIES

Assertion 6.

Let the arithmetic function $f(n), n = 1, 2, ....$ and the Dirichlet series $\sum_{n=1}^{\infty} \frac{f(n)}{n^s}$ - converge ($s \geq 1$) and the Dirichlet series $\sum_{n=1}^{\infty} \frac{f(n)}{n^{s-1}}$ - diverge, then:

$$\sum_{n \le x} \frac{f(n)}{n^{s-1}} = o(x) \qquad (3.1)$$

at $x \to \infty$.

Proof

If the Dirichlet series $\sum_{n=1}^{\infty} \frac{f(n)}{n^s}$ - converges, then $\sum_{n \le x} \frac{f(n)}{n^s} = c_1 + o(1)$ for the value $x \to \infty$.

Since the Dirichlet series $\sum_{n=1}^{\infty} \frac{f(n)}{n^{s-1}}$ - diverges, then taking into account (2.2) the following is true:

$$\sum_{n \le x} \frac{f(n)}{n^{s-1}} = \sum_{n \le x} (\frac{f(n)}{n^s})n = (c_1 + o(1))x - \int_1^x (c_1 + o(1))dt = c_1 x + o(x) - c_1 x + o(x) = o(x),$$

for the value $x \to \infty$, which corresponds to (3.1).

Let's look at an example for assertion 6.

It is known that the series $\sum_{n=1}^{\infty} \frac{\mu(n)}{n}$ - converges, and the series $\sum_{n=1}^{\infty} \mu(n)$ - diverges, therefore, based on assertion 6, is true:

$$\sum_{n \le x} \mu(n) = o(x), x \to 0, \qquad (3.2)$$

which is equivalent to the asymptotic law of prime numbers.

Asymptotic behavior (3.2) is well known. In other cases, similar asymptotics are not well known, so Assertion 7 [2] is useful for checking them.

Assertion 7

Let $|g(n)| \le 1, n = 1, ..., x$ is a real multiplicative arithmetic function. If the series $\sum_p \frac{1 - g(p)}{p} = \infty$ (diverges), then the asymptotic behavior of the summation function is equal to:

$$\sum_{n \le x} g(n) = o(x), x \to \infty. \qquad (3.3)$$

Let's look at a few more examples for assertion 6.

It is known that the series $\sum_{n=1}^{\infty} \frac{|\mu(n)|}{n^2} = \frac{\zeta(2)}{\zeta(4)}$, i.e. converges, but the series $\sum_{n=1}^{\infty} \frac{|\mu(n)|}{n}$ - diverges, therefore, based on assertion 6, it is true:

$$\sum_{n \leq x} \frac{|\mu(n)|}{n} = o(x), x \to \infty. \tag{3.4}$$

Let us check the asymptotic (3.4) based on assertion 7.

An arithmetic function $g(n) = \frac{|\mu(n)|}{n} \leq 1$ is multiplicative. The value $g(p) = \frac{|\mu(p)|}{p} = \frac{1}{p}$, therefore, $\sum_p \frac{1-g(p)}{p} = \sum_p \frac{1-1/p}{p} = \sum_p \frac{1}{p} - \sum_p \frac{1}{p^2}$ - diverges and based on assertion 7 (3.4) is also satisfied.

It is known that the series $\sum_{n=1}^{\infty} \frac{\sigma(n)}{n^3} = \zeta(2)\zeta(3)$, i.e. converges, but the series $\sum_{n=1}^{\infty} \frac{\sigma(n)}{n^2}$ - diverges, therefore, based on assertion 6, it is true:

$$\sum_{n \leq x} \frac{\sigma(n)}{n^2} = o(x), x \to \infty. \tag{3.5}$$

Let us give an example of using assertion 6 not for a multiplicative function. Let us show that the series $\sum_{n=1}^{\infty} \frac{\Omega(n)}{n^2}$ - converges:

$$\sum_{n=1}^{\infty} \frac{\Omega(n)}{n^2} \sim \sum_{n=1}^{\infty} \frac{\ln \ln n}{n^2} < \sum_{n=1}^{\infty} \frac{1}{n^{1+\xi}} \text{ - converges,}$$

where $\xi > 0$.

Now we show that the series $\sum_{n=1}^{\infty} \frac{\Omega(n)}{n}$ - diverges:

$$\sum_{n=1}^{\infty} \frac{\Omega(n)}{n} > \sum_{n=1}^{\infty} \frac{1}{n} \text{ - diverges.}$$

Therefore, based on assertion 6:

$$\sum_{n \leq x} \frac{\Omega(n)}{n} = o(x), x \to \infty. \tag{3.6}$$

Let's check (3.6):

$$\sum_{n \leq x} \frac{\Omega(n)}{n} \sim \sum_{n \leq x} \frac{\ln \ln n}{n} \sim \int_2^x \frac{\ln \ln t \, dt}{t} = \ln x (\ln \ln x - 1) ,$$

which corresponds to (3.6).

Corollary 8

Let there is an arithmetic function $f(n), n = 1, 2, ....$ and the Dirichlet series $\sum_{n=1}^{\infty} \frac{f(n)}{n^s}$ - converges $(s \geq 1)$, and the Dirichlet series $\sum_{n=1}^{\infty} \frac{f(n)}{n^{s-1}}$ - diverges, then the limit of the average value of the arithmetic function $\frac{f(n)}{n^{s-1}}$ is equal to:

$$\lim_{x \to \infty} E[\frac{f(n)}{n^{s-1}}, x] = \lim_{x \to \infty} \frac{1}{x} \sum_{n \leq x} \frac{f(n)}{n^{s-1}} = 0 . \tag{3.7}$$

Proof

The conditions of assertion 6 are satisfied, therefore:

$$\sum_{n \leq x} \frac{f(n)}{n^{s-1}} = o(x), x \to \infty .$$

Hence:

$$\lim_{x \to \infty} E[\frac{f(n)}{n^{s-1}}, x] = \lim_{x \to \infty} \frac{1}{x} \sum_{n \leq x} \frac{f(n)}{n^{s-1}} = \lim_{x \to \infty} \frac{o(x)}{x} = 0 ,$$

which corresponds to (3.7).

Let's look at an example for corollary 8.

The series $\sum_{n=1}^{\infty} \frac{\mu(n)}{n}$ - converges, and the series $\sum_{n=1}^{\infty} \mu(n)$ - diverges, therefore, based on corollary 8:

$$\lim_{x \to \infty} E[\mu(n), x] = \lim_{x \to \infty} \frac{1}{x} \sum_{n \leq x} \mu(n) = 0 . \tag{3.8}$$

In conclusion, I would like to note that many nonlinear asymptotic behavior of summation functions given in the work are known, but the methods for obtaining them in the existing literature are different and quite complex. Therefore, the goal of this work was to obtain a unified approach that would simplify obtaining these estimates.

4. CONCLUSION AND SUGGESTIONS FOR FURTHER WORK

The next article will continue to study the asymptotic behavior of some arithmetic functions.

5. ACKNOWLEDGEMENTS

Thanks to everyone who has contributed to the discussion of this paper. The author is very pleased that Professor Gérald Tenenbaum showed interest in this article and read it in full before publication in the Archive. I am grateful to everyone who expressed their suggestions and comments in the course of this work.